\documentclass[preprint,12pt]{elsarticleHF}
\usepackage{ngerman}
\usepackage{amsmath, amssymb, amsthm} %mehr Mathematikfunktionen
\usepackage{mathtools,mathdots} %enthält coloneqq
\numberwithin{equation}{section}
\usepackage{hyperref}	%Links erstellen
\usepackage[english]{cleveref} %cleveres Referenzieren, erstellt Worte wie "`Abschnitt"' automatisch
\usepackage{dsfont} %enthält die doppelt gestrichene 1 via \mathds{1}, ``double stroke''

\newtheorem{thm}{Theorem}
\newtheorem{cor}[thm]{Corollary}
\newtheorem{rem}[thm]{Remark}%[section]
\newtheorem{ex}[thm]{Example}

\theoremstyle{definition}

\DeclareMathOperator{\diag}{diag}

\begin{document}
\begin{frontmatter}

\title{On the singular value decomposition of (skew-)involutory and (skew-)coninvolutory matrices}

\author[label1]{Heike Fa\ss bender\corref{cor1}}
\address[label1]{Institut \emph{Computational Mathematics}/ AG Numerik, TU Braunschweig, Universitätsplatz 2, 38106 Braunschweig, Germany}
\cortext[cor1]{Corresponding author, Email: h.fassbender@tu-braunschweig.de }
\author[label2]{Martin Halwa\ss }
\address[label2]{Kampstraße 7, 17121 Loitz, Germany}

\begin{abstract}
The singular values $\sigma >1$ of an $n \times n$ involutory matrix $A$ appear in pairs
$(\sigma, \frac{1}{\sigma}).$ Their left and right singular vectors
are closely connected. The case of singular values $\sigma = 1$ is discussed in detail.
These singular values may appear in pairs $(1,1)$ with
closely connected left and right singular vectors or by themselves.
The link between the left and right singular vectors is used to reformulate the singular value decomposition (SVD) of an involutory matrix as
an eigendecomposition. This displays an interesting relation between the singular values of an involutory matrix and its eigenvalues.
Similar observations hold for the SVD, the singular values and the coneigenvalues of (skew-)coninvolutory matrices.
\end{abstract}

\end{frontmatter}

%%%%%%%%%%%%%%%%%%%%%%%%%%%%%%%%%%%%%%%%%%%%%%%%%%%%%%%%%%%%%%%
\section{Introduction}
Inspired by the work \cite{HouC63} on the singular values of involutory matrices some more insight into
the singular value decomposition (SVD) of involutory matrices is derived.
For any matrix $A\in \mathbb{C}^{n \times n}$ there exists a singular value decomposition (SVD), that is, a decomposition of the form
\begin{equation}\label{eq:svd} A = U \Sigma V^H,\end{equation}
where $U, V\in \mathbb{C}^{n\times n}$ are unitary matrices and $\Sigma \in \mathbb{R}^{n\times n}$  is a diagonal matrix
with non-negative real numbers on the diagonal. The diagonal entries $\sigma_{i}$ of $\Sigma$ are the singular values of $A.$
Usually, they are ordered such that $\sigma_1 \geq \sigma_2 \geq \cdots \geq \sigma_n \geq 0.$ The number of nonzero singular
values of $A$  is the same as the rank of $A.$ Thus, a nonsingular matrix $A \in \mathbb{C}^{n\times n}$ has  $n$ positive singular values.
The columns $u_j, j = 1, \ldots,n$  of $U$  and the
columns $v_j, j = 1, \ldots,n$  of $V$  are the left singular vectors and right singular vectors of $A,$ respectively.
 From \eqref{eq:svd} we
have $Av_j = \sigma_j u_j, j = 1, \ldots, n.$ Any triplet $(u,v,\sigma)$ with $Av = \sigma u$ is called a singular triplet of $A.$
In case $A \in \mathbb{R}^{n \times n},$ $U$ and $V$ can be chosen to be real and orthogonal.

While the singular values are
unique, in general, the singular vectors are not. The nonuniqueness of the singular vectors mainly depends on the multiplicities of the
singular values. For simplicity, assume that $A \in \mathbb{C}^{n \times n}$ is nonsingular. Let $s_1 > s_2 > \cdots > s_k >0$ denote the
distinct singular values of $A$ with respective multiplicities $\theta_1, \theta_2, \ldots, \theta_k \geq 1, \sum_{j=1}^k \theta_j = n.$
Let $A = U\Sigma V^H$ be a given singular value decomposition with $\Sigma = \diag(s_1I_{\theta_1}, s_2I_{\theta_2}, \ldots, s_kI_{\theta_k}).$
Here $I_j$ denotes as usual the $j \times j$ identity matrix. Then
\begin{equation}\label{nonuniqe_svd}
\widehat{U} = U \left[ W_1 \oplus W_2 \oplus \cdots \oplus W_k\right] \quad \text{and} \quad
\widehat{V} = V \left[ W_1 \oplus W_2 \oplus \cdots \oplus W_k\right]
\end{equation}
with unitary matrices $W_j \in \mathbb{C}^{\theta_j \times \theta_j}$ yield another SVD of $A,$ $A = \widehat{U}\Sigma \widehat{V}^H.$
This describes all possible SVDs of $A,$ see, e.g., \cite[Theorem 3.1.1']{HorJ94} or \cite[Theorem 4.28]{Ste98}. In case $\theta_j = 1,$ the corresponding left and right singular vector
are unique up to multiplication with some $e^{\imath \alpha_j}, \alpha_j \in \mathbb{R},$ where $\imath = \sqrt{-1}.$
For more information on the SVD see, e.g. \cite{GolV13,HorJ94,Ste98}.

A matrix $A\in \mathbb{C}^{n \times n}$ with $A^2=I_n,$ or equivalently,
$A = A^{-1}$ is called an involutory matrix.
Thus, for any involutory matrix $A,$
$A$ and its inverse $A^{-1}$ have the same SVD and hence, the same (positive) singular values.  Let $A = U\Sigma V^H$ be
the usual SVD of $A$ with the diagonal of $\Sigma$ ordered by magnitude, $\sigma_1 \geq \sigma_2 \geq \ldots \geq \sigma_n.$
Noting that an SVD of $A^{-1}$ is given by  $ A^{-1} = V \Sigma^{-1} U^H$ and that the diagonal elements of $\Sigma^{-1}$ are ordered by magnitude,
$\frac{1}{\sigma_n} \geq \frac{1}{\sigma_{n-1}} \geq \ldots \geq \frac{1}{\sigma_1},$ we conclude that
$\sigma_1 = \frac{1}{\sigma_n}, \sigma_2 = \frac{1}{\sigma_{n-1}}, \ldots, \sigma_n = \frac{1}{\sigma_1}.$
That is, the singular values of an involutory matrix $A$ are either $1$ or pairs $(\sigma, \frac{1}{\sigma}),$ where $\sigma >1.$
This has already been observed in \cite{HouC63} (see Theorem \ref{theo:householder} for a quote of the findings).
Here we will describe the SVD $A=U\Sigma V^H$ for involutory matrices in detail. In particular, we will note a close relation between the left and right singular vectors of pairs $(\sigma, \frac{1}{\sigma})$ of singular values:  if $(u,v,\sigma)$ is a singular triplet, then so is $(v,u,\frac{1}{\sigma}).$
We will observe that some of the  singular values $\sigma = 1$ may also appear in such pairs. That is, if $(u,v,1)$ is a singular triplet, then so is $(v,u,1).$ Other singular values $\sigma = 1$ may appear as a single singular triplet,
that is $(u, \pm u,1).$ These observations allow to express $U$ as $U =VT$ for  a real elementary orthogonal matrix $T.$
With this
the SVD of $A$ reads $A = VT\Sigma V^H.$ Taking a closer look at the real matrix $T\Sigma$ we will see that $T\Sigma$
is an involutory matrix just as $A.$ All relevant information concerning the singular values and the eigenvalues of $A$ can be read off of $T\Sigma.$
Some of these findings also follow from \cite[Theorem 7.2]{HorS09}, see Section \ref{sec:invol}.

As any skew-involutory matrix $B\in \mathbb{C}^{n \times n}, B^2 = -I_n$ can be expressed as $B = \imath A$ with an involutory matrix $A
\in \mathbb{C}^{n \times n},$ the results for involutory matrices can be transferred easily to skew-involutory ones.

We will also consider the SVD of coninvolutory matrices, that is of matrices $A \in \mathbb{C}^{ n \times n}$ which satisfy
$ A\overline{A} = I_n,$ see, e.g., \cite{HorJ13}. As involutory matrices, coninvolutory are nonsingular and have $n$ positive singular values.
Moreover, the singular values appear in pairs $(\sigma, \frac{1}{\sigma})$ or are $1,$ see \cite{HorJ94}.  Similar to the case of involutory matrices,
we can give a relation between the matrices $U$ and $V$ in the SVD $U\Sigma V$ of $A$ in the form
$U = \overline{V}T.$
$T\Sigma$ is a complex coninvolutory matrix consimilar\footnote{Two matrices $A$ and
$B\in \mathbb{C}^{n \times n}$  are said to be consimilar  if there exists a nonsingular matrix $S \in \mathbb{C}^{n \times n}$
such that $A=SB\overline{S}^{-1},$ \cite{HorJ13}.} to $A$ and consimilar to the identity. Similar observations have been given in \cite{HorJ13}.
Some of our findings also follow from \cite[Theorem 7.1]{HorS09}, see Section \ref{sec:coninvol}.

Skew-coninvolutory matrices $A \in \mathbb{C}^{2n \times 2n},$ that is, $A\overline{A}=-I_{2n},$ have been studied in \cite{AbaMP07}.
Here we will briefly state how with our approach findings on the SVD of skew-coninvolutory matrices given in  \cite{AbaMP07} can easily be rediscovered.

Our goal is to round off the picture drawn in the literature about the singular value decomposition of the four classes of matrices considered here.
In particular we would like to make visible the relation between the singular vectors belonging to reciprocal pairs of singular values in the form of the
matrix $T.$

The SVD of involutory matrices is treated at length in Section \ref{sec:invol}.
In Section \ref{sec:skewinvol} we make immediate use of the results in Section \ref{sec:invol} in order to discuss the SVD of skew-involutory matrices.
Section \ref{sec:coninvol} deals with coninvolutory matrices. Finally, the SVD of skew-coninvolutory matrices is discussed in  Section \ref{sec:skewconinvol}.

%%%%%%%%%%%%%%%%%%%%%%%%%%%%%%%%%%%%%%%%%%%%%%%%%%%%%%%%%%
\section{Involutory matrices}\label{sec:invol}
%%%%%%%%%%%%%%%%%%%%%%%%%%%%%%%%%%%%%%%%%%%%%%%%%%%%%%%%%%
Let $A\in \mathbb{C}^{n \times n}$  be involutory, that is, $A^2=I_n$ holds. Thus, $A = A^{-1}.$
Symmetric permutation matrices and Hermitian unitary matrices are simple examples of involutory matrices.
Nontrivial examples of involutory matrices can be found in \cite[Page 165, 166, 170]{Hig08}.
The spectrum of any involutory matrix can have at most two
elements, $\lambda(A) \subset \{1, -1\}$ as  from $Ax = \lambda x$ for $x\in \mathbb{C}^n\backslash{0}$ it follows that
$A^2x = \lambda Ax$ which gives $x = \lambda^2x.$

 The SVD of involutory matrices has already been considered in \cite{HouC63}. In particular, the following theorem is given.
\begin{thm}%[\cite{HouC63}]
\label{theo:householder}
Let $A\in \mathbb{C}^{n \times n}$ be an involutory matrix. Then $B_1 =\frac{1}{2}( I_n-A)$ and $B_2 = \frac{1}{2}(I_n+A)$ are idempotent,
that is, $B_i^2=B_i, i = 1,2.$
Assume that $A$ has $r \leq \frac{n}{2}$ eigenvalues $-1$ and $n-r$ eigenvalues $+1$ (or $r \leq \frac{n}{2}$ eigen\-values $+1$ and $n-r$ eigenvalues $-1$). Then $B_1$ ($B_2$) is of rank $r$ and
can be decomposed as $B_1 = QW^H$ ($B_2 = QW^H$) where $Q, W \in \mathbb{C}^{n \times r}$ have $r$ linear independent columns.
Moreover, $W^HW-I$ is at least positive semidefinite.
There are $n-2r$ singular values of $A$ equal to $1,$ $r$ singular values  which are equal to the eigenvalues
of the matrix $(W^HW)^\frac{1}{2}+(W^HW-I_n)^\frac{1}{2}$ and $r$ singular values which are the reciprocal of these.
\end{thm}

Thus, in \cite{HouC63} it was noted that the singular values of an involutory matrix $A$ may be $1$ or may appear in pairs
$(\sigma, \frac{1}{\sigma}).$ Moreover, the minimal number of singular values $1$ is given by $2n-r$ where $r \leq \frac{n}{2}$
denotes the number of eigenvalues $-1$ of $A$ or the number of eigenvalues $+1$ of $A$, whichever is smallest.
Further, in \cite{HouC63} it is said ''that the singular values $\neq 1$ of the involutory matrix $A$ are the roots of''
$(W^HW)^\frac{1}{2}+(W^HW-I_n)^\frac{1}{2}$ (and $(W^HW)^\frac{1}{2}-(W^HW-I_n)^\frac{1}{2}$).
This does not imply that
the matrix $(W^HW)^\frac{1}{2}+(W^HW-I_n)^\frac{1}{2}$ does not have eigenvalues equal to $+1.$
This can be illustrated by $A = \operatorname{diag}(+1, -1, -1)$ with one eigenvalue $+1$ and two eigenvalues $-1.$ Hence, in
Theorem \ref{theo:householder}, $r = 1$ and there is (at least) one  singular value $1,$ as $n-2r = 1.$ For the other two singular values,  the matrix
$B_2 = \frac{1}{2}(I_n+A) = QW^T$ with $Q=[1, 0, 0]^T$ and $W^T=[1, 0, 0]$ and the $1 \times 1$ matrix
$(W^TW)^\frac{1}{2}+(W^TW-I_n)^\frac{1}{2} = 1$ needs to be considered. Obviously, $W$ has
the eigenvalue $+1$ which, according to Theorem \ref{theo:householder}  is a singular value of $A.$ Moreover, its reciprocal has to be a singular value.
Hence, $A$ has three singular values $1,$  two appearing as a pair $(\sigma,\frac{1}{\sigma})= (1,1)$  and a 'single' one.

In the following discussion we will see that apart from the pairing of the singular values, there is even more structure in the SVD of an involutory matrix
by taking a closer look at $U$ and $V.$ This will highlight the difference between the two types of singular values $1$ (pairs of singular values $(1,1)$
and single singular values $1$).

Let $(u,v,\sigma)$ be a singular triplet of $A$, that is, let  $Av = \sigma u$ hold. This is equivalent to $v = \sigma A^{-1}u = \sigma A u$ and further to $Au = \frac{1}{\sigma}v.$ Thus, the
singular triplet  $(u,v,\sigma)$ of $A$ is always accompanied by the singular triplet  $(v,u,\frac{1}{\sigma})$ of $A.$
 In case $\sigma = 1,$ this may collapse into a single triplet if $u=v$ or $u=-v.$
 This allows for three cases:
\begin{enumerate}
\item  In case $u=v,$ we have a single singular value $1$ with the  triplet $(u,u,1).$ In this case, it follows immediately that $A$ has an eigenvalue $1,$ as we have $Au = u.$
\item  In case $u=-v,$ we have a single singular value $1$ with the  triplet $(u,-u,1)$ (or the triplet $(-u,u,1),$). In this case, it follows immediately that $A$ has an eigenvalue $-1,$ as we have $Au = -u.$
\item In case $u \neq \pm v,$ we have a pair $(1,1)$ of singular values associated with the two triplets $(u,v,1)$ and $(v,u,1).$
\end{enumerate}
Therefore, singular values $\sigma > 1$ will appear in pairs $(\sigma, \frac{1}{\sigma}),$ while a singular value $\sigma =1$ may appear in a pair in the sense
that there are two triplets $(u,v,1)$ and $(v,u,1)$ or it may appear by itself as a triplet $(u,\pm u,1).$
It follows immediately, that an involutory matrix $A$ of odd size $n$ must have at least one singular value $1.$

Clearly, due to the nonuniqueness of the SVD, no every SVD of $A$ has to display the pairing of the singular values identified above.
But every SVD of $A$ can be easily modified so that this can be read off.  Let us consider a small example first.
\begin{ex}\label{ex}
Consider the involutory, Hermitian and unitary matrix
\[ A = \begin{bmatrix} 1 & 0&0&0\\ 0 &0&1&0\\ 0&1&0&0\\ 0&0&0&1\end{bmatrix} \in \mathbb{R}^{4 \times 4}.\]
One possible SVD of $A$ is given by $A = \widetilde{U} \Sigma \widetilde{V}^T$ with
\[ \widetilde{U} = \begin{bmatrix}
\frac{1}{\sqrt{2}} & \frac{\sqrt{3}}{2\sqrt{2}}  & \frac{1}{2\sqrt{2}}&0\\
0& -\frac{1}{2} & \frac{\sqrt{3}}{2} & 0\\
- \frac{1}{\sqrt{2}}  & \frac{\sqrt{3}}{2\sqrt{2}} & \frac{1}{2\sqrt{2}} & 0\\
0 & 0 & 0 & 1
\end{bmatrix}, \qquad\Sigma = I_4,\qquad
\widetilde{V} = \begin{bmatrix}
\frac{1}{\sqrt{2}} & \frac{\sqrt{3}}{2\sqrt{2}} & \frac{1}{2\sqrt{2}}&0\\
- \frac{1}{\sqrt{2}}  & \frac{\sqrt{3}}{2\sqrt{2}} & \frac{1}{2\sqrt{2}} & 0\\
0& -\frac{1}{2} & \frac{\sqrt{3}}{2} & 0\\
0 & 0 & 0 & 1
\end{bmatrix}.
\]
We have  $A\widetilde{v}_1 = \widetilde{u}_1, A\widetilde{v}_2 = \widetilde{u}_2,$ $A\widetilde{v}_3 =\widetilde{u}_3,$
 and $A\widetilde{v}_4  = Ae_4 =e_4= \widetilde{u}_4.$
It can be seen immediately that there is  the singular triplet $(1,e_4,e_4).$
Thus, there needs to be at least one other singular triplet of this type, but it can not be seen straightaway whether the other two
singular values $1$ correspond to a pair of singular values with connected singular vectors or not.
 Please note that
$\widetilde{U} = UW,  \widetilde{V} =VW$ holds for
\begin{align*}
U = A=\begin{bmatrix} e_1, e_3, e_2, e_4 \end{bmatrix}, \quad
V = I_4 = \begin{bmatrix} e_1, e_2, e_3, e_4 \end{bmatrix}, \quad
W = \widetilde{V}.
\end{align*}
Thus $A = U \Sigma V^H$ is another SVD of $A$ with
\begin{align*}
 Ae_1 = e_1,\quad
 Ae_2 = e_3,\quad
 Ae_3 = e_2,\quad
Ae_4 = e_4.
\end{align*}
Hence, there is one singular value which is part of the two related triplets $(e_2,e_3,1)$ and $(e_3,e_2,1)$ and two singular
values which have no partner as their triplets are $(e_1,e_1,1)$ and $(e_4,e_4,1).$
 In a similar way, any other SVD of $A$ can be modified in order to display the structure of the singular values and vectors.
\end{ex}

Assume that
an SVD $A = U \Sigma V^H$ is given where the singular values are ordered such that
$\sigma_1 \geq \sigma_2 \geq \cdots \geq \sigma_n >0.$ Let us first assume that $\sigma_1 > \sigma_2.$
Then $\sigma_{n-1}= \frac{1}{\sigma_2}> \sigma_n = \frac{1}{\sigma_1}.$ Moreover, $Av_1 = \sigma_1u_1$ and
$Av_n=\sigma_nu_n = \frac{1}{\sigma_1}u_n,$ or, equivalently, $Av_1 = \sigma_1u_1$ and $Au_n = \sigma_1v_n.$
Due to  the essential uniqueness of singular vectors
of singular values with multiplicity $1,$ it follows that $u_n = v_1e^{\imath \alpha_1}$ and $v_n = u_1e^{\imath \alpha_1}$
for some $\alpha_1.$ Modifying $U$ and $V$ to
\begin{align*}
\tilde U & := \left[ u_1~\mid~u_2~~\cdots~~u_{n-1}~\mid~v_1\right]\\
\tilde V & := \left[ v_1~\mid~v_2~~\cdots~~v_{n-1}~\mid~u_1\right]
\end{align*}
yields unitary matrices $\tilde U$ and $\tilde V$ such that $A = \tilde U\Sigma \tilde V^H$ is a valid SVD displaying the
pairing for the singular value $\sigma_1.$ In this fashion all singular values $>1$ of multiplicity $1$ can be treated.

Next let us assume that $\sigma_1= \sigma_2 = \cdots = \sigma_\ell > \sigma_{\ell+1}.$
Then $\sigma_{n-\ell}= \frac{1}{\sigma_{\ell+1}}>\sigma_{n-\ell+1} = \cdots= \sigma_n = \frac{1}{\sigma_1}.$
Due to our observation concerning the singular vectors
of pairs of singular values $(\sigma, \frac{1}{\sigma})$ and the essential uniqueness of the SVD \eqref{nonuniqe_svd},
it follows that
\begin{align*}
\left[ u_{n-\ell+1}~~\cdots~~u_{n-1}~~u_n\right] &= \left[ v_\ell~~\cdots~~v_2~~v_1\right]W_1\\
\left[ v_{n-\ell+1}~~\cdots~~v_{n-1}~~v_n\right] &= \left[ u_\ell~~\cdots~~u_2~~u_1\right]W_1
\end{align*}
for some unitary $\ell \times \ell$ matrix $W_1.$ Modifying $U$ and $V$ to
\begin{align*}
\tilde U & := \left[ u_1~~u_2~~\cdots~~u_\ell~\mid~u_{\ell+1}~~\cdots ~~u_{n-\ell}~\mid~v_\ell~~\cdots~~v_2~~v_1\right]\\
\tilde V & := \left[ v_1~~v_2~~\cdots~~v_\ell~\mid~v_{\ell+1}~~\cdots ~~v_{n-\ell}~\mid~u_\ell~~\cdots~~u_2~~u_1\right]
\end{align*}
yields unitary matrices $\tilde U$ and $\tilde V$ such that $A = \tilde U\Sigma \tilde V^H$ is a valid SVD displaying the
pairing for the singular value $\sigma_1= \sigma_2 = \cdots = \sigma_\ell.$ In this fashion all singular values $>1$ of multiplicity $>1$
can be treated.

Finally we need to consider the singular values $1.$ Similar as before, we can modify the columns of $U$ and $V$
so that the relation between the singular vectors becomes apparent.

This gives rise to the following theorem.
\begin{thm}[SVD of an involutory matrix]\label{theo1}
Let $A\in \mathbb{C}^{n \times n}$ be involutory.
Assume that $A$ has $\nu$ singular values $\sigma_1\geq \cdots \geq \sigma_\nu >1.$
These singular values appear in pairs $(\sigma_j,\frac{1}{\sigma_j})$ associated with the singular triplets
 $(u_j,v_j,\sigma_j)$  and $(v_j,u_j,\frac{1}{\sigma_j}), j = 1, \ldots, \nu.$
Assume further that $A$ has $\mu$ singular values $1$ which
appear in pairs $(1,1)$ associated with the singular triplets
 $(\tilde{u}_j,\tilde{v}_j,1)$  and $(\tilde{v}_j,\tilde{u}_j,1), j = 1, \ldots, \mu.$
Finally, assume that $A$ has $k=n-2\nu-2\mu$ single singular values $1$  associated with the singular triplet
 $(\hat{u}_j,\pm\hat{u}_j,1), j = 1, \ldots, k.$

Thus the SVD of $A$ is given by
\[ \Sigma = \begin{bmatrix} S &&&\\ & I_\mu & \\ &&I_\delta & \\&&& S^{-1}&\\ &&&&I_\mu \\ &&&&&I_\eta\end{bmatrix}
 \]
with $S = \diag(\sigma_1, \ldots, \sigma_\nu)$ and
\begin{align*}
U &= \left[ u_1, \ldots, u_\nu, \tilde{u}_1, \ldots, \tilde{u}_\mu, \hat{u}_1, \ldots, \hat{u}_\delta, \right.\\
&\qquad\quad \left. v_1, \ldots, v_\nu, \tilde{v}_1, \ldots, \tilde{v}_\mu, \hat{u}_{\delta+1}, \ldots, \hat{u}_{\eta+\delta}\right],\\
V &= \left[ v_1, \ldots, v_\nu, \tilde{v}_1, \ldots, \tilde{v}_\mu, \pm\hat{u}_1, \ldots, \pm\hat{u}_\delta, \right.\\
&\qquad \quad \left. u_1, \ldots, u_\nu, \tilde{u}_1, \ldots, \tilde{u}_\mu, \pm\hat{u}_{\delta+1}, \ldots, \pm\hat{u}_{\eta+\delta}\right],
\end{align*}
where $ \nu+\mu+\eta=m$ and
\[
\delta = \left\{ \begin{array}{ll}
 \eta & \text{if } n =2m \\
 \eta+1 & \text{if } n=2m+1
\end{array} \right. .
\]
In particular, the signs do not need to be equal for all $\hat{u}_j, j =1, \ldots, \delta+\eta.$
\end{thm}

For $U$ and $V$ from Theorem \ref{theo1} we have
\begin{equation}\label{eq:T}
 U = V \begin{bmatrix}
0 & 0 & 0 & I_\nu & 0 & 0\\
0 & 0 & 0 & 0 & I_\mu & 0 \\
0 & 0 & D_{\delta} & 0 & 0 & 0\\
I_\nu & 0 & 0 & 0 & 0 & 0\\
0 & I_\mu & 0 & 0 & 0 & 0\\
0 & 0 & 0 & 0 & 0 &  E_\eta
\end{bmatrix}= V \begin{bmatrix}
 0 & 0 & I_{\nu+\mu}  & 0\\
 0 & D_{\delta}  & 0 & 0\\
I_{\nu+\mu}  & 0 & 0  & 0\\
 0 & 0 & 0 &  E_\eta
\end{bmatrix} = VT
\end{equation}
where $D_\delta\in \mathbb{R}^{\delta \times \delta}, E_\eta\in \mathbb{R}^{\eta \times \eta}$ denote diagonal matrices with $\pm 1$
 on the diagonal. The particular choice depends on the sign choice in the sequence $\pm\hat{u}_j$ in $V$ in Theorem \ref{theo1}.
Clearly, $D_\delta$ and $E_\eta$ as well as $T$
are involutory.

Thus we have
\begin{equation}\label{eq:TSigma}
 A = U \Sigma V^H
= VT\Sigma V^H=
V  \begin{bmatrix}
0 & 0 & 0 & S^{-1} & 0 & 0\\
0 & 0 & 0 & 0 & I_\mu & 0 \\
0 & 0 & D_{\delta} & 0 & 0 & 0\\
S & 0 & 0 & 0 & 0 & 0\\
0 & I_\mu & 0 & 0 & 0 & 0\\
0 & 0 & 0 & 0 & 0 &  E_\eta
\end{bmatrix} V^H.
\end{equation}
In other words, $A$ is unitarily similar to the real involutory matrix $T\Sigma.$ This canonical form is the most
condensed involutory matrix unitarily similar to $A.$ All relevant information concerning the singular values and
the eigenvalues of $A$ can be read off of $T\Sigma.$

Making use of the fact that all diagonal elements of $S$ are positive, we can rewrite
$S$ as $S = S^{\frac{1}{2}} S^{\frac{1}{2}}$ with $ S^{\frac{1}{2}} =
\diag(\sqrt{\sigma_1}, \sqrt{\sigma_2},\ldots, \sqrt{\sigma_\nu}).$ Thus
\begin{align*}
A &=V
\begin{bmatrix}
S^{-\frac{1}{2}} & 0 & 0 & 0 & 0 & 0\\
0 & I_\mu & 0 & 0 & 0 & 0 \\
0 & 0 & I_{\delta} & 0 & 0 & 0\\
0 & 0 & 0 & S^{\frac{1}{2}} & 0 & 0\\
0 & 0 & 0 & 0 & I_\mu & 0\\
0 & 0 & 0 & 0 & 0 &  I_\eta
\end{bmatrix}
T
\begin{bmatrix}
S^{\frac{1}{2}} & 0 & 0 & 0 & 0 & 0\\
0 & I_\mu & 0 & 0 & 0 & 0 \\
0 & 0 & I_{\delta} & 0 & 0 & 0\\
0 & 0& 0 & S^{-\frac{1}{2}} & 0  & 0\\
0 & 0 & 0 & 0 & I_\mu & 0\\
0 & 0 & 0 & 0 & 0 & I_\eta
\end{bmatrix}
V^H\\
&=: ZTZ^{-1}.
\end{align*}
Hence, $A$ and $T$ are similar matrices and their eigenvalues are identical. Taking a closer look at $T,$ we immediately see that
$T$
 is similar to
\begin{equation}\label{eq:Diag}
PTP^T=
\begin{bmatrix}  -I_{\nu+\mu} &0&0& 0\\
0&I_{\nu+\mu}&0&0\\
0 & 0&D_\delta & 0\\
0&0&0&E_\eta\end{bmatrix}
\end{equation}
with the orthogonal matrix
\[
P = \begin{bmatrix} \frac{1}{\sqrt{2}}I_{\nu+\mu} & 0&-\frac{1}{\sqrt{2}}I_{\nu+\mu}&0\\
\frac{1}{\sqrt{2}}I_{\nu+\mu}&0&\frac{1}{\sqrt{2}}I_{\nu+\mu}&0\\
0 & I_\delta & 0 & 0\\
0&0&0&I_\eta\end{bmatrix}.
\]
Assume that there are $\eta_1$ positive and $\eta_2$ negative signs in the sequence $\pm\hat{u}_j$ in $V,$ $j = 1, \ldots, \delta+\eta=\eta_1+\eta_2.$  Then $\diag(D_\delta, E_\eta)$ is similar to $\diag(I_{\eta_1}, -I_{\eta_2}).$
 It is straightforward to see that
\[
\det(T-\lambda I_n) = (\lambda+1)^{\nu+\mu}(\lambda-1)^{\nu+\mu}(\lambda-1)^{\eta_1} (\lambda+1)^{\eta_2}.
\]
Each pair of singular triplets $(u,v,\sigma)$ and $(v,u,\frac{1}{\sigma})$ (including those with $\sigma =1$)
corresponds to a pair of eigenvalues $(+1,-1).$
A single singular triplet $(u,\pm u, 1)$ corresponds to an eigenvalue $+1$ or $-1$ depending on the sign in $(u,\pm u, 1).$
Our findings are summarized in the following corollary.
\begin{cor}[Canonical Form, Eigendecomposition]\label{cor1}
Let $A\in \mathbb{C}^{n \times n}$ be involutory. Let $A = U \Sigma V^H$ be the SVD of $A.$
Assume that $A$ has $\nu$ singular values $>1,$ $\mu$ singular values $1$ which
appear in pairs $(1,1),$  $\eta_1$ single singular values $1$  associated with the singular triplet
 $(\hat{u}_j,\hat{u}_j,1)$ and $\eta_2$ single singular values $1$  associated with the singular triplet
 $(\hat{u}_j,-\hat{u}_j,1).$ Here $\eta_1+\eta_2=\delta+\eta$ for $\delta, \eta$ as in Theorem \ref{theo1}.
Then $A$ is unitarily similar to the real involutory matrix $T\Sigma$ as in \eqref{eq:TSigma} and diagonalizable to
$\diag(-I_{\nu+\mu+\eta_2}, I_{\nu+\mu+\eta_1}),$ see \eqref{eq:Diag}.
\end{cor}

\begin{rem}
If $A$ is involutory, then $B = \frac{1}{2}(I\pm A)$ is idempotent, $B^2 = B.$ This has been used in \cite{HouC63}, see Theorem \ref{theo:householder}.
Using \eqref{eq:TSigma} we obtain for $B$ that $B =  \frac{1}{2}(I\pm VT\Sigma V^H) = \frac{1}{2}V(I\pm T\Sigma )V^H= \frac{1}{2}VT(T\pm \Sigma )V^H$ holds as $T$ is involutory.
We can easily construct an SVD of $T\pm \Sigma,$
\[
T \pm \Sigma %=
%\begin{bmatrix}
%0 & 0 & 0 & I_\nu & 0 & 0\\
%0 & 0 & 0 & 0 & I_\mu & 0 \\
%0 & 0 & D_{\delta} & 0 & 0 & 0\\
%I_\nu & 0 & 0 & 0 & 0 & 0\\
%0 & I_\mu & 0 & 0 & 0 & 0\\
%0 & 0 & 0 & 0 & 0 &  E_\eta
%\end{bmatrix}
%  \pm
%\begin{bmatrix} S &&&\\ & I_\mu & \\ &&I_\delta & \\&&& S^{-1}&\\ &&&&I_\mu \\ &&&&&I_\eta\end{bmatrix}
= \begin{bmatrix}
\pm S & 0 & 0 & I_\nu & 0 & 0\\
0 & \pm I_\mu & 0 & 0 & I_\mu & 0 \\
0 & 0 & D_{\delta} \pm I_\delta& 0 & 0 & 0\\
I_\nu & 0 & 0 & \pm S^{-1}  & 0 & 0\\
0 & I_\mu & 0 & 0 & \pm I_\mu & 0\\
0 & 0 & 0 & 0 & 0 &  E_\eta \pm I_\eta
\end{bmatrix}.
\]
First permute $T\pm \Sigma$ with
\[ P_1 =\begin{bmatrix}
I_\nu &0 & 0 & 0 & 0 & 0\\
0 & 0 & I_\mu & 0 & 0 & 0\\
0& 0& 0 & 0 & I_\delta & 0\\
0 & I_\nu & 0 & 0 &0&0\\
0 & 0 & 0 & I_\mu & 0 & 0\\
0& 0 &0 & 0 & 0 & I_\eta
\end{bmatrix}
\]
such that
\[
P_1^T(T\pm \Sigma)P_1 =
\begin{bmatrix}
\pm S & I_\nu & 0& 0 & 0 & 0\\
I_\nu & \pm S^{-1} &0&0 & 0 & 0\\
0 & 0 & \pm I_\mu & I_\mu & 0 & 0\\
0 & 0 & I_\mu  & \pm I_\mu & 0 & 0\\
0 & 0 & 0 & 0 & D_\delta \pm I_\delta & 0\\
0 & 0 & 0 & 0 & 0 & E_\eta \pm I_\eta
\end{bmatrix}.
\]
Next permute the  blocks
$\left[\begin{smallmatrix} \pm S & I_\nu\\ I_\nu & \pm S^{-1}\end{smallmatrix}\right]$
and $\left[\begin{smallmatrix} \pm I_\mu & I_\mu\\ I_\mu & \pm I_\mu\end{smallmatrix}\right]$
to block diagonal form $\widetilde{S} = \diag\left( \left[\begin{smallmatrix} \pm \sigma_1 & 1\\ 1 & \pm \sigma_1^{-1}\end{smallmatrix}\right], \ldots,
\left[\begin{smallmatrix} \pm \sigma_\nu & 1\\ 1 & \pm \sigma_\nu^{-1}\end{smallmatrix}\right] \right)$
and $\widetilde{I} = \diag\left(  \left[\begin{smallmatrix} \pm 1 & 1\\ 1& \pm 1\end{smallmatrix}\right], \ldots, \left[\begin{smallmatrix} \pm 1 & 1\\ 1& \pm 1\end{smallmatrix}\right]  \right) .$ Let $P_2$ be the corresponding permutation matrix such that
\[
P_2^TP_1^T (T\pm \Sigma)P_1 P_2=
\begin{bmatrix}
\widetilde{S} & 0& 0 & 0 \\
0 & \widetilde{I}&0 & 0 \\
 0 & 0 & D_\delta \pm I_\delta & 0\\
 0 & 0 & 0 & E_\eta \pm I_\eta
\end{bmatrix}
\]
is block diagonal with $1 \times 1$ and $2 \times 2$ diagonal blocks.
The $2 \times 2$ blocks are real symmetric and can be diagonalized by an orthogonal similarity transformation
\[
\begin{bmatrix}\pm \sigma& 1\\ 1 & \pm \sigma^{-1}\end{bmatrix} =
\begin{bmatrix}c & s\\  -s& c\end{bmatrix}
\begin{bmatrix}\pm(\sigma+ \sigma^{-1})&0 \\ 0& 0\end{bmatrix}
\begin{bmatrix}c& -s\\  s& c\end{bmatrix} ,
\]
with
$c = \sqrt{\frac{\sigma}{\sigma+\sigma^{-1}}}$ and $s = \frac{c}{\sigma}.$
Let $X$ be the orthogonal matrix which diagonalizes $P_2^TP_1^T (T\pm \Sigma)P_1 P_2,$
\[
X^TP_2^TP_1^T (T\pm \Sigma)P_1 P_2X =
\begin{bmatrix}
\pm \widehat{S} & 0& 0 & 0 \\
0 & \pm \widehat{I}&0 & 0 \\
 0 & 0 & D_\delta \pm I_\delta & 0\\
 0 & 0 & 0 & E_\eta \pm I_\eta
\end{bmatrix}
\]
with $\widehat{S} = \diag ( \sigma_1+\sigma_1^{-1},0, \ldots,\sigma_\nu+\sigma_\nu^{-1},0)$ and $\widehat{I} = \diag(2,0, \ldots,2,0 ).$
This gives an SVD of $ T\pm \Sigma$ and thus of  $B.$ In case $B = \frac{1}{2}(I+A)$ has been chosen, we have
\[
 T + \Sigma =  P_1P_2X \diag(\widehat{S}, \widehat{I}, D_\delta + I_\delta , E_\eta+ I_\eta) X^TP_2^TP_1^T,
\]
and 
\[
B = \underbrace{VTP_1P_2X}_{unitary} 
\underbrace{\diag(\frac{1}{2}\widehat{S}, \frac{1}{2}\widehat{I}, \frac{1}{2}(D_\delta+ I_\delta) , \frac{1}{2}(E_\eta + I_\eta))}_{nonnegative~ diagonal} \underbrace{X^TP_2^TP_1^TV^H}_{unitary}.
\]
In case $B = \frac{1}{2}(I-A)$ has been chosen, we need to take care of the minus sign in front of $\widehat{S}$ and $\widehat{I}$
\begin{align*}
 T - \Sigma &=  P_1P_2X \diag(-\widehat{S}, -\widehat{I}, D_\delta - I_\delta , E_\eta - I_\eta) X^TP_2^TP_1^T\\
&= P_1P_2X \diag(\widehat{S}, \widehat{I}, D_\delta - I_\delta , E_\eta - I_\eta) Y^TP_2^TP_1^T 
\end{align*}
with $Y = -X.$ The SVD of $B$ is immediate.

Hence, if $A$ has pairs of singular values $(\sigma, \sigma^{-1})$ and $(1,1),$ then $B$ has pairs of singular values
$(\frac{\sigma+\sigma^{-1}}{2},0)$ and $(1,0).$ Moreover, $D_\delta \pm I_\delta$ and $E_\eta \pm I_\eta$ give singular values $1$ or $0.$
\end{rem}

Before we turn our attention to the skew-involutory case, we would like to point out that most of our observations given in this section also
follow  from  \cite[Theorem 7.2]{HorS09}. For the ease of the reader, this theorem is stated next.
\begin{thm} \label{theo:hornsergeichuk}
Let $A\in \mathbb{C}^{n \times n}.$ If $A^2$ is normal, then $A$ is unitarily $H$-congruent\footnote{Two matrices $A$ and
$B \in \mathbb{C}^{n \times n}$ are said to be $H$-congruent if there exists a nonsingular matrix $S \in \mathbb{C}^{n \times n}$
such that $A = SBS^H.$ Thus, two unitarily $H$-congruent matrices are unitarily similar.} to a direct sum of blocks, each of which is
\begin{equation}\label{eq_xx}
 [\lambda] \quad \text{ or } \quad \tau \begin{bmatrix} 0 & 1\\ \mu & 0\end{bmatrix}, \quad \tau \in \mathbb{R}, \lambda, \mu \in \mathbb{C}, \tau > 0, \text{ and } |\mu|< 1.
\end{equation}
This direct sum is uniquely determined by $A$, up to permutation of its blocks. Conversely, if $A$
is unitarily $H$-congruent to a direct sum of blocks of the form \eqref{eq_xx}, then $A^2$ is normal.
\end{thm}
In case $A^2 = I,$ Theorem \ref{theo:hornsergeichuk} gives for the $1 \times 1$ blocks that $\lambda = \pm 1$ has to hold.
For the $2 \times 2$ blocks it follows  $\tau^2\mu = 1.$ Thus, $\mu \in (0,1)$ and $\tau\mu = \frac{1}{\tau}.$
Let $U \in \mathbb{C}^{n \times n}$
be the unitary matrix which transforms the involutory matrix $A$ as described in Theorem  \ref{theo:hornsergeichuk} ;
\begin{equation}\label{eq_xxx}
U^HAU =\big( \bigoplus_{i=1}^{n_1} [\lambda_i]\big) \oplus \left( \bigoplus_{j=1}^{n_2}
\tau_j \begin{bmatrix} 0 & 1\\ \mu_j & 0\end{bmatrix}\right),
\end{equation}
with $ \tau_j, \mu_j \in \mathbb{R}, \lambda_j \in \{+1,-1\}, \tau_j > 0, \mu_j \in (0,1)$ and $ \tau_j\mu_j = \frac{1}{\tau_j}.$
Clearly, $\tau_j \neq 1$ as $\mu_j \neq 1.$

The unitary $H$-congruence of $A$ as in \eqref{eq_xxx} can be modified into an SVD.
For any $1 \times 1$ block $[\lambda_i]$, the corresponding singular value is $1.$ At the same time,
any $1 \times 1$ block $[\lambda_i]$ represent an eigenvalue of $A.$ An eigenvector $u_i$ corresponding to $\lambda_i$
can be read off of $U; Au_i = \lambda_i u_i.$ This eigenvector $u_i$ will serve as the corresponding right singular vector.
The left singular vector will be chosen as $u_i$ in case the $ \lambda_i =1$ and as $-u_i$ in case $\lambda_i = -1.$
The SVD of a  $2 \times 2$ block is given by
\[
\tau \begin{bmatrix} 0 & 1\\ \mu  & 0\end{bmatrix} = I_2 \begin{bmatrix} \tau  & 0\\0 & \tau  \mu\end{bmatrix}
\begin{bmatrix} 0 & 1\\1 & 0\end{bmatrix} = W \Sigma Z^H.
\]
Thus, as $\tau\mu = \frac{1}{\tau}$ and $\mu \in (0,1)$ holds, the singular values $\tau >1$ appear in pairs $(\tau, \frac{1}{\tau}).$
There are columns $v, w$ from $U$ such that
\[
A[v ~~w] = [v ~~w] \tau
\begin{bmatrix} 0 & 1\\ \mu & 0\end{bmatrix} = [v ~~w] \Sigma Z^H.
\]
 This gives
\[
A[w ~~v] = [v ~~w]\begin{bmatrix} \tau & 0\\ 0 & \frac{1}{\tau}\end{bmatrix}.
\]
Hence, singular values $\tau >1$ appear in pairs $(\tau,\frac{1}{\tau})$
and are associated with the
singular triplets $(v, w,\tau)$ and $(w,v,\frac{1}{\tau}).$
The fact that singular values $1$ can also appear in pairs $(1,1)$ with singular triplets in the form $(u,v,1)$ and $(v,u,1)$
does not follow from Theorem \ref{theo:hornsergeichuk}.

The unitary $H$-congruence of $A$ can also be modified into an eigendecomposition.
The $1 \times1$ blocks represent eigenvalues of $A,$ a corresponding eigenvector can be read off of $U.$  The eigenvalues of
the $2 \times 2$ blocks are $+1$ and $-1.$ Each $2 \times 2$ block can be diagonalized by a unitary matrix.
Thus, Theorem  \ref{theo:hornsergeichuk}  yields that any involutory matrix is unitarily
diagonalizable to $\diag(-I_{n_{11}+n_2},+I_{n_{12}+n_2}).$  Here, it is assumed
%The eigendecomposition of a $2 \times 2$ block is given by
%\[
%\tau_j \begin{bmatrix} 0 & 1\\ \mu_j & 0\end{bmatrix} = XX \begin{bmatrix} 1 & 0\\ 0 & -1\end{bmatrix} YY
%\]
 that there are $n_2$ $2 \times 2$ blocks and $n_1 = n-2n_2$ $1 \times 1$ blocks with $n_{11}$ blocks $[-1]$ and $n_{12}$ blocks $[+1].$
Comparing this result to our one in Corollary \ref{cor1} we see that $\nu = n_2, \mu+\eta_2 = n_{11}$ and $\mu+\eta_1 = n_{12}.$

%%%%%%%%%%%%%%%%%%%%%%%%%%%%%%%%%%%%%%%%%%%%%%%
\section{Skew-involutory matrices}\label{sec:skewinvol}
%%%%%%%%%%%%%%%%%%%%%%%%%%%%%%%%%%%%%%%%%%%%%%%
Any skew-involutory matrix $A \in \mathbb{C}^{n \times n}$ can be expressed as $A = \imath C$ with an involutory matrix
$C\in \mathbb{C}^{n \times n}.$  Thus we can immediately make use of the
results from Section \ref{sec:invol}. As for the spectrum of an involutory matrix $C$ we have $\lambda(C) \subset \{1, -1\},$ it follows that
$\lambda(A) \subset \{\imath,-\imath\}$ holds. Moreover, if the singular value decomposition of $C$ is given by $C = U \Sigma V^H$ with
$U, \Sigma, V$ as in Theorem \ref{theo1}, then $A = U \Sigma (\imath V^H)$ is an SVD of $A.$

In particular, it holds that  any singular value $\sigma >1$ of $A$ appears as a pair $(\sigma, \frac{1}{\sigma}).$ The
singular triplet  $(u,v,\sigma)$ of $A$ is always accompanied by the singular triplet  $(-v,u,\frac{1}{\sigma})$ of $A.$
 In case $\sigma = 1,$ this may collapse into a single triplet if $u=\imath v$ or $u=-\imath v.$ Thus, the SVD of $A$ can be given as in Theorem \ref{theo1} where $V$ is modified such that
\begin{align*}
V &= \left[ v_1, \ldots, v_\nu, \tilde{v}_1, \ldots, \tilde{v}_\mu, \pm \imath \hat{u}_1, \ldots, \pm \imath \hat{u}_{\delta},\right.\\
&\qquad\quad \left. -u_1, \ldots, -u_\nu, -\tilde{u}_1, \ldots, -\tilde{u}_\mu, \pm\imath  \hat{u}_{\delta+1}, \ldots, \pm\imath\hat{u}_{\eta+\delta}\right].
\end{align*}
Similar to before, $U$ and $V$ are closely connected
\[
 U =
V \begin{bmatrix}
 0 & 0 & -I_{\nu+\mu}  & 0\\
 0 & \imath D_{\delta}  & 0 & 0\\
I_{\nu+\mu}  & 0 & 0  & 0\\
 0 & 0 & 0 &  \imath E_\eta
\end{bmatrix} = VT.
\]
Hence, we have
\[
A = U \Sigma V^H
= VT\Sigma V^H=
V  \begin{bmatrix}
0 & 0 & 0 & -S^{-1} & 0 & 0\\
0 & 0 & 0 & 0 & -I_\mu & 0 \\
0 & 0 & \imath D_{\delta} & 0 & 0 & 0\\
S & 0 & 0 & 0 & 0 & 0\\
0 & I_\mu & 0 & 0 & 0 & 0\\
0 & 0 & 0 & 0 & 0 &  \imath E_\eta
\end{bmatrix} V^H.
\]
In other words, $A$ is unitarily similar to the  complex skew-involutory matrix $T\Sigma$ which reveals all relevant information about
the singular values and the eigenvalues of $A.$
Moreover, $A$ is diagonalizable  to $\diag( -\imath I_{\nu+\mu+\eta_2},\imath I_{\nu+\mu+\eta_1})$
with $\eta_1, \eta_2$ as in Corollary \ref{cor1}.

Please note, that Theorem \ref{theo:hornsergeichuk} holds for skew-involutory matrices. Similar comments as those given at the end of
Section \ref{sec:invol} hold here.

 %%%%%%%%%%%%%%%%%%%%%%%%%%%%%%%%%%%%%%%%%%%%%%%
\section{Coninvolutory matrices}\label{sec:coninvol}
 %%%%%%%%%%%%%%%%%%%%%%%%%%%%%%%%%%%%%%%%%%%%%%%
For any coninvolutory matrix $A \in \mathbb{C}^{n \times n}$ we have $A^{-1} = \overline{A}$ as $A\overline{A}=I_n.$
Any coninvolutory matrix can be expressed as $A = e^{\imath R}$ for $R \in \mathbb{R}^{n\times n},$ see, e.g., \cite{HorJ94}.
Any real coninvolutory matrix is also involutory. Since  $A^{-1} = \overline{A},$ the singular values of $A$ are either $1$ or pairs $\sigma, \frac{1}{\sigma}.$
Moreover, any coninvolutory matrix is condiagonalisable\footnote{ $A  \in \mathbb{C}^{n \times n}$ is said to be condiagonalizable if there exists a nonsingular matrix $S \in \mathbb{C}^{n \times n}$ such that $D=S^{-1}A\overline{S}$ is diagonal \cite{HorJ13}.}, see, e.g., \cite[Chapter 4.6]{HorJ13}.

Let $(u,v,\sigma)$ be  singular triplet of a coninvolutory matrix $A$, that is, let  $Av = \sigma u$ hold. This is equivalent to $\overline{v} = \sigma \overline{A}^{-1}\overline{u} = \sigma A \overline{u}$ and further to $A\overline{u} = \frac{1}{\sigma}\overline{v}.$ Thus,  the
singular triplet  $(u,v,\sigma)$ of $A$ is always accompanied by the singular triplet  $(\overline{v},\overline{u},\frac{1}{\sigma})$ of $A.$
 In case $\sigma = 1,$ this may collapse into a single triplet if $v=e^{\imath \alpha} \overline{u}$  for a
real scalar $\alpha \in [0,2\pi]$ (as $-e^{\imath \alpha} = e^{\imath (\alpha+\pi) },$ there is no need to consider $v = -e^{\imath \alpha}\overline{u}$).
The case $\sigma = 1,$ $v=e^{\imath \alpha}\overline{u}$ implies that $A$ has a coneigenvalue $e^{-\imath \alpha}$ as $ A\overline{u}=e^{-\imath \alpha}u$ holds\footnote{A nonzero vector $x \in \mathbb{C}^n$ such that $A\overline{x} = \lambda x$ for some $\lambda\in \mathbb{C}$ is said to be an coneigenvector of $A;$
the scalar $\lambda$ is an coneigenvalue of $A,$ \cite{HorJ13}.}.
It follows immediately, that coninvolutory matrix $A$ of odd size $n$ must have a singular value $1.$

This gives rise the following theorem.
\begin{thm}\label{theo2}
Let $A\in \mathbb{C}^{n \times n}$ be coninvolutory.
Assume that $A$ has $\nu$ singular values $\sigma_1\geq \cdots \geq \sigma_\nu >1.$
These singular values appear in pairs $(\sigma_j,\frac{1}{\sigma_j})$ associated with the singular triplets
 $(u_j,v_j,\sigma_j)$  and $(\overline{v}_j,\overline{u}_j,\frac{1}{\sigma_j}), j = 1, \ldots, \nu.$
Assume further that $A$ has $\mu$ singular values $1$ which
appear in pairs $(1,1)$ associated with the singular triplets
 $({u}_{\nu+j},{v}_{\nu+j},1)$  and $(\overline{{v}}_{\nu+j},\overline{{u}}_{\nu+j},1), j = 1, \ldots, \mu.$
Then $A$ has $k=n-2\nu-2\mu=\delta + \eta$ single singular values $1$  associated with the singular triplet
 $({u}_{\nu+\mu+j}, e^{\imath \alpha_j}\overline{{u}}_{\nu+\mu+j},1), j = 1, \ldots, k.$

Thus the SVD of $A$ is given by
\[ \Sigma = \begin{bmatrix} S &&&\\ & I_\mu & \\ &&I_\delta & \\&&& S^{-1}&\\ &&&&I_\mu \\ &&&&&I_\eta\end{bmatrix}
\]
with $S = \diag(\sigma_1, \ldots, \sigma_\nu)$ and
\begin{align*}
U &= \left[ u_1, \ldots, u_\nu, {u}_{\nu+1}, \ldots, {u}_{\nu+\mu}, {u}_{\nu+\mu+1}, \ldots, {u}_{\nu+\mu+\delta}, \right.\\
&\qquad\quad \left. \overline{v}_1, \ldots, \overline{v}_\nu, \overline{{v}}_{\nu+1}, \ldots, \overline{{v}}_{\nu+\mu}, {u}_{\nu+\mu+\delta+1}, \ldots, {u}_{\nu+\mu+\delta +\eta}\right],\\
V &= \left[ v_1, \ldots, v_\nu, {v}_{\nu+1}, \ldots, {v}_{\nu+\mu},  e^{\imath \alpha_1}\overline{{u}}_{\nu+\mu+1}, \ldots,  e^{\imath \alpha_\delta}\overline{{u}}_{\nu+\mu+\delta},\right.\\
&\qquad\quad \left. \overline{u}_1, \ldots, \overline{u}_\nu, \overline{{u}}_{\nu+1}, \ldots, \overline{{u}}_{\nu+\mu},  e^{\imath \alpha_{\delta+1}}\overline{{u}}_{\nu+\mu+\delta+1}, \ldots,  e^{\imath \alpha_{\delta+\eta}}\overline{{u}}_{\nu+\mu+\delta+{\eta}}\right]
\end{align*}
where $ \nu+\mu+\eta=m$ and
\[
\delta = \left\{ \begin{array}{ll}
 \eta & \text{if } n =2m \\
 \eta+1 & \text{if } n=2m+1
\end{array} \right. .
\]
\end{thm}

For $U$ and $V$ from Theorem \ref{theo2} we have $U = \overline{V}T$
\begin{equation}\label{eq:Ttheo2}
 U =
\overline{V} \begin{bmatrix}
 0 & 0 & I_{\nu+\mu}  & 0\\
 0 & D_{\delta}  & 0 & 0\\
I_{\nu+\mu}  & 0 & 0  & 0\\
 0 & 0 & 0 &  E_\eta
\end{bmatrix} = \overline{V}T
\end{equation}
with  the unitary and coninvolutory diagonal matrices
\begin{align*}
D_{\delta} &= \diag( e^{\imath \alpha_1}, \ldots,  e^{\imath \alpha_\delta}),\\
E_{\eta} &= \diag( e^{\imath \alpha_{\delta+1}}, \ldots,  e^{\imath \alpha_{\delta+\eta}}).
\end{align*}
Thus we have
$ A = U \Sigma V^H = \overline{V}T\Sigma V^H$ with
\begin{equation}\label{eq:TSigma_coninvol}
T\Sigma =
\begin{bmatrix}
0 & 0 & 0 & S^{-1} & 0 & 0\\
0 & 0 & 0 & 0 & I_\mu & 0 \\
0 & 0 & D_{\delta} & 0 & 0 & 0\\
S & 0 & 0 & 0 & 0 & 0\\
0 & I_\mu & 0 & 0 & 0 & 0\\
0 & 0 & 0 & 0 & 0 &  E_\eta
\end{bmatrix}.
\end{equation}
In other words, $A$ is unitarily consimilar to the complex coninvolutory matrix $T\Sigma.$
A similar statement is given in \cite[Exercise 4.6P27]{HorJ13} (just write
$D_{\delta}$ and $ E_\eta$ as a product of their square roots and move the square roots into $\overline{V}$ and $V^H$).

As in Section \ref{sec:invol} we obtain
\begin{align*}
A &=\overline{V}
\begin{bmatrix}
S^{-\frac{1}{2}} & 0 & 0 & 0 & 0 & 0\\
0 & I_\mu & 0 & 0 & 0 & 0 \\
0 & 0 & I_{\delta} & 0 & 0 & 0\\
0 & 0 & 0 & S^{\frac{1}{2}} & 0 & 0\\
0 & 0 & 0 & 0 & I_\mu & 0\\
0 & 0 & 0 & 0 & 0 &  I_\eta
\end{bmatrix}
T
\begin{bmatrix}
S^{\frac{1}{2}} & 0 & 0 & 0 & 0 & 0\\
0 & I_\mu & 0 & 0 & 0 & 0 \\
0 & 0 & I_{\delta} & 0 & 0 & 0\\
0 & 0& 0 & S^{-\frac{1}{2}} & 0  & 0\\
0 & 0 & 0 & 0 & I_\mu & 0\\
0 & 0 & 0 & 0 & 0 &  I_\eta
\end{bmatrix}
V^H\\
&=: \overline{Z}TZ^{-1}.
\end{align*}
Hence, $A$ and $T$ are consimilar matrices and their coneigenvalues are identical. Taking a closer look at $T,$ we immediately see that
$T$
 is unitarily consimilar to the identity
\[
PTP^T=
\begin{bmatrix}  I_{\nu+\mu} &0 &0& 0\\
0& I_{\delta}&0&0\\
0 & 0&I_{\nu+\mu} & 0\\
0&0&0&I_\eta\end{bmatrix}
\]
with the unitary matrix
\[
P = \begin{bmatrix} \frac{1}{\sqrt{2}}I_{\nu+\mu} & 0&  \frac{1}{\sqrt{2}}I_{\nu+\mu}&0\\
0 & D_\delta^{-\frac{1}{2}}& 0 & 0\\
 -\frac{\imath}{\sqrt{2}} I_{\nu+\mu}&0&  \frac{\imath}{\sqrt{2}}I_{\nu+\mu}&0\\
0&0&0&E_\eta^{-\frac{1}{2}}\end{bmatrix}.
\]
Thus, all coneigenvalues of a coninvolutory matrix are $+1.$ This has already been observed in \cite[Theorem 4.6.9]{HorJ13}.

The following corollary summarizes our findings.
\begin{cor}[Canonical Form, Coneigendecomposition]\label{cor1_coninvol}
Let $A\in \mathbb{C}^{n \times n}$ be coninvolutory.
Then $A$ is unitarily consimilar to the  coninvolutory matrix $T\Sigma$ as in \eqref{eq:TSigma_coninvol} and consimilar to the identity.
\end{cor}

Before we turn our attention to the skew-coninvolutory case, we would like to point out that most of our observations given in this section also
follow easily from  \cite[Theorem 7.1]{HorS09}. For the ease of the reader, this theorem is stated next.
\begin{thm} \label{theo:hornsergeichuk2}
Let $A\in \mathbb{C}^{n \times n}.$ If $\overline{A}A$ is normal, then $A$ is unitarily congruent to a direct sum of blocks, each of which is
\begin{equation}\label{eq_yy}
 [\lambda] \quad \text{ or } \quad \tau \begin{bmatrix} 0 & 1\\ \mu & 0\end{bmatrix}, \quad \lambda,\tau \in \mathbb{R}, \lambda \geq 0, \tau > 0,
\mu \in \mathbb{C}, \text{ and } \mu\neq 1.
\end{equation}
This direct sum is uniquely determined by $A$ up to permutation of its blocks and replacement of any nonzero parameter $\mu$ by $\mu^{-1}$ with a corresponding replacement of $\tau$ by $\tau|\mu|.$ Conversely, if $A$
is unitarily congruent to a direct sum of blocks of the form \eqref{eq_yy}, then $\overline{A}A$ is normal.
\end{thm}
In case $\overline{A}A = I,$ Theorem \ref{theo:hornsergeichuk2} gives for the $1 \times 1$ blocks that $\lambda = \pm 1$ has to hold.
For the $2 \times 2$ blocks it follows  $\tau^2\mu = 1$ and $\tau^2\overline{\mu}=1.$ Thus, as $\mu \neq 1$ we have $\tau \neq 1.$
Analogous to the discussion at the end of Section \ref{sec:invol} part of our findings, in particular the pairing $(\tau, \frac{1}{\tau})$ of the singular values $\tau >1$
and the relation of the corresponding singular vectors, follows from this; see also \cite[Corollary 8.4]{HorS09}. The fact, that
singular values $1$ may also appear in pairs $(1,1)$ with related singular triplets does not follow from Theorem \ref{theo:hornsergeichuk2}.
 %%%%%%%%%%%%%%%%%%%%%%%%%%%%%%%%%%%%%%%%%%%%%%%
\section{Skew-coninvolutory matrices}\label{sec:skewconinvol}
%%%%%%%%%%%%%%%%%%%%%%%%%%%%%%%%%%%%%%%%%%%%%%%
For any skew-coninvolutory matrix $A \in \mathbb{C}^{m \times m}$ we have $A^{-1} = -\overline{A}$ as $A\overline{A}=-I_m.$
Skew-coninvolutory matrices exist only for even $m,$ as $\det A\overline{A}$ is nonnegative for any
$A \in \mathbb{C}^{m \times m}.$
Properties and canonical forms of skew-coninvolutory matrices have been analyzed in detail in \cite{AbaMP07}.

From $A = U\Sigma V^H$ we see that $A^{-1} = V\Sigma^{-1}U^H = -\overline{U}
\Sigma V^T.$ Thus, the  singular values appear in pairs $\sigma, \frac{1}{\sigma}$ (see \cite[Proposition 5]{AbaMP07}) and  the
singular triplet  $(u,v,\sigma)$ of $A$ is always accompanied by the singular triplet  $(-\overline{v},\overline{u},\frac{1}{\sigma})$ of $A.$
There is no need to consider
singular values $\sigma = 1$ separately, as their singular vectors do not satisfy any additional condition.
This gives rise to the following theorem.
\begin{thm}\label{theo3}
Let $A\in \mathbb{C}^{2n \times 2n}$ be skew-coninvolutory.
 $A$ has $n$ singular values $\sigma_1\geq \cdots \geq \sigma_n \geq1.$
These singular values appear in pairs $(\sigma_j,\frac{1}{\sigma_j})$ associated with the singular triplets
 $(u_j,v_j,\sigma_j)$  and $(-\overline{v}_j,\overline{u}_j,\frac{1}{\sigma_j}), j = 1, \ldots, n.$

Thus the SVD of $A$ is given by
\[ \Sigma = \begin{bmatrix} S & 0\\ 0 & S^{-1}\end{bmatrix} \]
with $S = \diag(\sigma_1, \ldots, \sigma_n)$ and
\begin{align*}
U &= \begin{bmatrix} u_1, \ldots, u_n, -\overline{v}_1, \ldots, -\overline{v}_n \end{bmatrix},\\
V &= \begin{bmatrix} v_1, \ldots, v_n,  \overline{u}_1,  \ldots, \overline{{u}}_n \end{bmatrix}.
\end{align*}
\end{thm}

For $U$ and $V$ from Theorem \ref{theo3} we have $U = -\overline{V}J_n$ with $J_n = \left[\begin{smallmatrix} 0 & I_n \\ -I_n & 0\end{smallmatrix}\right].$
Thus we have
$ A = U \Sigma V^H = -\overline{V}J_n\Sigma V^H.$ Any skew-coninvolutory matrix is unitarily consimilar to
the elementary skew-coninvolutory matrix $-J_n\Sigma.$
Observe
\begin{equation}\label{eq:TSigma_skewconinvol}
J_n\Sigma =
\begin{bmatrix}
0 & S^{-1}\\ -S & 0
\end{bmatrix} =
\begin{bmatrix}
S^{-\frac{1}{2}}&0\\
0& S^{\frac{1}{2}}
\end{bmatrix}
J_n
\begin{bmatrix}
S^{\frac{1}{2}}&0\\
0&S^{-\frac{1}{2}}
\end{bmatrix},
\end{equation}
(see \cite[Theorem 12]{AbaMP07}).
Let $Z = V\left[\begin{smallmatrix} S^{-\frac{1}{2}} &0\\ 0&S^{\frac{1}{2}}  \end{smallmatrix}\right].$
Thus $A = -\overline{Z}J_n Z^{-1}$
 and $-J_n$ are consimilar matrices (see \cite[Theorem 13]{AbaMP07}).

In summary, we have the following corollary.
\begin{cor}[Canonical Form, Consimilarity]\label{cor1_skewconinvol}
Let $A\in \mathbb{C}^{2n \times 2n}$ be skew-con\-involutory.
Then $A$ is unitarily consimilar to the  skew-coninvolutory matrix $-J_n\Sigma$ as in \eqref{eq:TSigma_skewconinvol} and consimilar to $-J_n.$
\end{cor}

Please note, that Theorem \ref{theo:hornsergeichuk2} holds for skew-coninvolutory matrices. Similar comments as those given at the end of
Section \ref{sec:coninvol} hold here, see also \cite[Theorem 8.3]{HorS09}.
%%%%%%%%%%%%%%%%%%%%%%%%%%%%%%%%%%%%%%%%%%%%%%%
\section{Concluding remarks}\label{sec:conclud}
We have described the SVD of (skew-)involutory and (skew-)coninvolutory matrices in detail. In order to do so we made use of the fact that
for any matrix in one of the four classes of matrices
the singular values $\sigma >1$ appear in pairs $(\sigma, \frac{1}{\sigma}),$ while singular values $\sigma = 1$
may appear in pairs or by themselves. As the singular vectors of pairs of singular values are closely related, the SVD
reveals all relevant information also about the eigenvalues and eigenvectors. Some of our findings are new, some are rediscoveries of known results.
%%%%%%%%%%%%%%%%%%%%%%%%%%%%%%%%%%%%%%%%%%%%%%%

%\bibliographystyle{elsarticle-harv}
%\bibliographystyle{abbrv}
%\bibliographystyle{plain}
%\bibliography{mybib}

\begin{thebibliography}{1}

\bibitem{AbaMP07}
Ma. Nerissa~M. Abara, Dennis~I. Merino, and Agnes~T. Paras.
\newblock Skew-coninvolutory matrices.
\newblock {\em Linear Algebra Appl.}, 426(2-3):540--557, 2007.

\bibitem{GolV13}
Gene~H. Golub and Charles~F. Van~Loan.
\newblock {\em Matrix computations}.
\newblock Johns Hopkins Studies in the Mathematical Sciences. Johns Hopkins
  University Press, Baltimore, MD, fourth edition, 2013.

\bibitem{Hig08}
Nicholas~J. Higham.
\newblock {\em Functions of matrices}.
\newblock Society for Industrial and Applied Mathematics (SIAM), Philadelphia,
  PA, 2008.
\newblock Theory and computation.

\bibitem{HorJ94}
Roger~A. {Horn} and Charles~R. {Johnson}.
\newblock {\em {Topics in matrix analysis. 1st pbk with corr.}}
\newblock Cambridge: Cambridge University Press, 1994.

\bibitem{HorJ13}
Roger~A. {Horn} and Charles~R. {Johnson}.
\newblock {\em {Matrix analysis. 2nd ed.}}
\newblock Cambridge: Cambridge University Press, 2013.

\bibitem{HorS09}
Roger~A. {Horn} and Vladimir~V. {Sergeichuk}.
\newblock Canonical forms for unitary congruence and *congruence.
\newblock {\em Linear Multilinear Algebra}, 57(8):777--815, 2009.

\bibitem{HouC63}
Alston~S. Householder and John~A. Carpenter.
\newblock The singular values of involutory and of idempotent matrices.
\newblock {\em Numer. Math.}, 5:234--237, 1963.

\bibitem{Ste98}
Gilbert~W. {Stewart}.
\newblock Matrix algorithms. {V}ol. {I},  {Basic decompositions}.
 \newblock Society for Industrial and Applied Mathematics, Philadelphia,
              PA, {1998}.

\end{thebibliography}

\end{document}